\documentclass[11pt]{article}
\usepackage{amscd}
\usepackage{amsfonts}
\usepackage{amsmath}
\usepackage{amssymb}
\usepackage{amsthm}
\usepackage{bbm}
\usepackage{CJK}
\usepackage{fancyhdr}
\usepackage{graphicx}
\usepackage{hyperref}
\usepackage{indentfirst}
\usepackage{latexsym}
\usepackage{mathrsfs}
\usepackage{xypic}
\usepackage{amsmath,amssymb,amsfonts,amsthm,fancyhdr}
\usepackage{epsfig,graphicx,picins,picinpar,subfigure}
\usepackage{pstricks}
\usepackage{amssymb}
\usepackage{xypic}
\usepackage[compress]{cite}

\newtheorem{theorem}{Theorem}[section]
\newtheorem{prop}[theorem]{Proposition}
\newtheorem{lemma}[theorem]{Lemma}
\newtheorem{coro}[theorem]{Corollary}
\newtheorem{prop-def}{Proposition-Definition}[section]

\theoremstyle{definition}
\newtheorem{definition}[theorem]{Definition}

\newtheorem{remark}[theorem]{Remark}
\newtheorem{exam}[theorem]{Example}

\usepackage[top=1in,bottom=1in,left=1.25in,right=1.25in]{geometry}
\def\<{\langle}
\def\>{\rangle}

\date{\today}
\begin{document}
\renewcommand{\baselinestretch}{1.2}
\renewcommand{\arraystretch}{1.0}
\title{\bf Embedding tensors on 3-Leibniz algebras and their  derived algebraic structures  and deformations }
\author{{\bf     Wen Teng, Shuangjian Guo }\\
{\small  School of Mathematics and Statistics, Guizhou University of Finance and Economics} \\
{\small  Guiyang  550025, P. R. of China}\\
 {E-mail: tengwen@mail.gufe.edu.cn (Wen Teng),}\\
  { shuangjianguo@126.com (Shuangjian Guo)}\\
 }
 \maketitle
\begin{center}
\begin{minipage}{13.cm}

{\bf Abstract:}
In this paper,   first we introduce the notions of   3-tri-Leibniz algebras and  embedding tensors on 3-Leibniz algebras.
We show that   an embedding tensor gives rise to a 3-tri-Leibniz algebra.
 Conversely,  a 3-tri-Leibniz algebra
gives rise to a 3-Leibniz algebra and a representation such that the quotient map is an embedding tensor.
  Furthermore,    any 3-tri-Leibniz algebra can be embedded into an averaging 3-Leibniz algebra.
Next, we introduce the notion  of 3-tri-Leibniz dialgebras and
 demonstrate that  homomorphic embedding tensors inherently induce   3-tri-Leibniz dialgebras.
Finally, we study the linear  deformations of embedding tensors  by defining the first  cohomology.
 \smallskip

 \smallskip

 \end{minipage}
 \end{center}
 \normalsize\vskip0.5cm

\section{Introduction}
\def\theequation{\arabic{section}. \arabic{equation}}
\setcounter{equation} {0}

Leibniz algebras can be viewed as a non-skew-symmetric counterpart to Lie algebras.
The class of $n$-Leibniz algebras, serving as a natural extension of Leibniz algebras to higher arities, was introduced in \cite{Casas}.
Furthermore, an $n$-Leibniz algebra, when combined with skew-symmetry, constitutes an $n$-Lie algebra \cite{Filippov}.
 In recent years, both domestic and international scholars have undertaken extensive research on the structure and properties of   $n$-Leibniz algebras.
 For instance,  in \cite{Albeverio},  the authors explored  the characteristics of Cartan subalgebras and normal
elements in $n$-Leibniz algebras. In  \cite{Casas16}, the authors employed Gr$\mathrm{\ddot{o}}$bner bases to formulate an algorithm for verifying the given multiplication table corresponding to an $n$-Leibniz algebra.
In  \cite{XU}, the authors studied non-abelian extensions of 3-Leibniz algebras through Maurer-Cartan elements.
See  \cite{Azcarraga}  for more details about $n$-Leibniz algebras.

Embedding tensors emerge in the exploration of supergravity theory and advanced gauge theories \cite{Nicolai}.
  From a   mathematical point of view, embedding tensors are also known as averaging operators \cite{Aguiar}.
  Recently, the mathematical theory of embedding tensors on Lie algebras, Lie$_\infty$-algebras,  3-Lie algebras, groups, Leibniz algebras,  mock-Lie
algebras and  Lie triple systems
 was studied in \cite{Sheng,Das24,Braiek,Caseiro,Hu,Teng25}. It is showed that an  embedding tensor on a 3-Lie algebra or a Lie triple system can give rise to
a 3-Leibniz algebra structure, see \cite{Hu,Teng25}.

Considering the  significance of 3-Leibniz algebras and embedding tensors,
 this paper seeks to present the concept of embedding tensors on 3-Leibniz algebras,
  and to explore the underlying algebraic structures and linear deformations associated with them.
More precisely, in Section  \ref{sec: 3-tri-Leibniz algebras},  we  introduce the notion of a 3-tri-Leibniz algebra and provide some examples.
In Section \ref{sec: Embedding tensors},  we  introduce the notions of  embedding tensors and averaging operators on
3-Leibniz algebras. We show that an embedding tensor or an  averaging operator induces a 3-tri-Leibniz algebra.
Conversely,  a 3-tri-Leibniz algebra gives rise to a 3-Leibniz algebra and a representation such that the quotient map is an embedding tensor.
Furthermore, we show that any
3-tri-Leibniz algebra can be embedded into an averaging 3-Leibniz algebra.
 In Section \ref{sec: 3-tri-Leibniz dialgebras},   we  introduce the notion of a 3-tri-Leibniz dialgebra and
show that  homomorphic embedding tensors induce naturally 3-tri-Leibniz dialgebras.
In Section \ref{sec: Deformations}, we define the first  cohomology to study the  linear deformations of
embedding tensors on 3-Leibniz algebras, following Gerstenhaber's approach \cite{Gerstenhaber}.

Throughout this paper, $\mathbb{K}$ denotes a field of characteristic zero. All the    vector spaces  and
   (multi)linear maps are taken over $\mathbb{K}$.

\section{ 3-tri-Leibniz algebras }\label{sec: 3-tri-Leibniz algebras}
\def\theequation{\arabic{section}.\arabic{equation}}
\setcounter{equation} {0}

This section introduces the notion of a 3-tri-Leibniz algebra and provides several examples.
To do this, we need to recall a representation of a 3-Leibniz algebra which allows the replicating of the operation into three parts with some compatibilities.

\begin{definition} \cite{Casas}
A 3-Leibniz algebra is a vector space $\mathfrak{g}$ endowed with a  trilinear map $[\cdot, \cdot, \cdot]_\mathfrak{g}:\mathfrak{g}\times\mathfrak{g}\times\mathfrak{g}\rightarrow\mathfrak{g}$,
which satisfies the following equation:
\begin{align}
 &[a, b, [x, y, z]_\mathfrak{g}]_\mathfrak{g}=[[a, b, x]_\mathfrak{g},y,z]_\mathfrak{g}+ [x,  [a, b, y]_\mathfrak{g},z]_\mathfrak{g}+ [x,y,[a, b, z]_\mathfrak{g}]_\mathfrak{g} \label{2.1}
\end{align}
for all $ x, y, z, a, b\in \mathfrak{g}$.
\end{definition}

Let  $(\mathfrak{g},[\cdot,\cdot,\cdot]_{\mathfrak{g}})$   be a 3-Leibniz algebra. A vector subspace of $\mathfrak{g}$
closed under the  operation $[\cdot,\cdot,\cdot]_{\mathfrak{g}}$ is
called a 3-Leibniz subalgebra of $\mathfrak{g}$. An ideal of $\mathfrak{g}$ is a vector subspace
$I$ which satisfies $[I, \mathfrak{g}, \mathfrak{g}]_{\mathfrak{g}}+[\mathfrak{g}, I, \mathfrak{g}]_{\mathfrak{g}}+[\mathfrak{g}, \mathfrak{g}, I]_{\mathfrak{g}}\subseteq I$.

\begin{definition}
Let $(\mathfrak{g},[\cdot,\cdot,\cdot]_{\mathfrak{g}})$ and $(\mathfrak{h},[\cdot,\cdot,\cdot]_{\mathfrak{h}})$ be two 3-Leibniz algebras. If
for any $ x, y, z\in \mathfrak{g}$, linear map $\varphi:\mathfrak{g}\rightarrow\mathfrak{h}$ satisfies
$$\varphi[x,y,z]_{\mathfrak{g}}=[\varphi(x),\varphi(y),\varphi(z)]_{\mathfrak{h}},$$
then $\varphi$ is said to be a 3-Leibniz algebra homomorphism.
\end{definition}

\begin{exam}
Any vector space with the trivial trilinear  bracket is a 3-Leibniz  algebra, called an abelian 3-Leibniz algebra.
\end{exam}

\begin{exam}
3-Lie algebras and Lie triple systems are 3-Leibniz algebras.
\end{exam}

\begin{exam}
$\mathbb{R}^4$ is a 3-Leibniz algebra with the bracket given by $[u, v, w]_{\mathbb{R}^4}=u\times v\times w$, where $u\times v\times w$ is the vector product of the vectors
$u,v,w\in\mathbb{R}^4$.
\end{exam}

\begin{exam}
Let  $(\mathfrak{g},[\cdot,\cdot])$  be a Leibniz algebra. Then $(\mathfrak{g},[\cdot,\cdot,\cdot]_{\mathfrak{g}})$ is a 3-Leibniz algebra, where $[x,y,z]_\mathfrak{g}=[[x,y],z]$  for all
$x,y,z\in \mathfrak{g}$. Conversely,  if $(\mathfrak{g},[\cdot,\cdot,\cdot]_{\mathfrak{g}})$ is a 3-Leibniz algebra, then  $\mathfrak{g}\otimes\mathfrak{g}$ with the operation
$$[\mathcal{X},\mathfrak{Y}]=[x_1,x_2,y_1]_{\mathfrak{g}}\otimes y_2+ y_1\otimes[x_1,x_2,y_2]_{\mathfrak{g}}~\text{for all}~\mathcal{X}=x_1\otimes x_2,\mathfrak{Y}=y_1\otimes y_2\in \mathfrak{g}\otimes\mathfrak{g}$$
is a Leibniz algebra,
\end{exam}

\begin{definition} \cite{Casas} A representation of the 3-Leibniz algebra $(\mathfrak{g},[\cdot,\cdot,\cdot]_{\mathfrak{g}})$ is a
vector space $V$ equipped with three trilinear maps:
\begin{align*}
&\rho^l:\mathfrak{g}\otimes \mathfrak{g}\otimes V\rightarrow V,~~~\rho^m:\mathfrak{g}\otimes V\otimes \mathfrak{g}\rightarrow V,~~~\rho^r:V\otimes \mathfrak{g}\otimes \mathfrak{g}\rightarrow V,
\end{align*}
such that for any  $x,y,z,a,b\in \mathfrak{g}$ and $u\in V$,
\begin{align}
&\rho^l(a,b,\rho^l(x,y,u))=\rho^l([a,b,x]_{\mathfrak{g}}, y, u)+ \rho^l(x,[a,b,y]_{\mathfrak{g}},u)+\rho^l(x,y,\rho^l(a,b,u)),\label{2.2}\\
 &\rho^l(a,b,\rho^m(x,u,z))= \rho^m([a,b,x]_{\mathfrak{g}}, u, z)+\rho^m(x,\rho^l(a,b,u),z)+\rho^m(x,u,[a,b,z]_{\mathfrak{g}}),\label{2.3}\\
  &\rho^l(a,b,\rho^r(u,y,z))=\rho^r(\rho^l(a,b,u), y,z)+ \rho^r(u,[a,b,y]_{\mathfrak{g}},z)+\rho^r(u,y,[a,b,z]_{\mathfrak{g}}),\label{2.4}\\
   &\rho^m(a,u,[x,y,z]_{\mathfrak{g}})=\rho^r(\rho^m(a,u,x), y,z)+ \rho^m(x,\rho^m(a,u,y),z)+\rho^l(x,y,\rho^m(a,u,z)),\label{2.5}\\
    &\rho^r(u,b,[x,y,z]_{\mathfrak{g}})=\rho^r(\rho^r(u,b,x), y,z)+ \rho^m(x,\rho^r(u,b,y),z)+\rho^l(x,y,\rho^r(u,b,z)).\label{2.6}
\end{align}
\end{definition}

\begin{exam}
Given a 3-Leibniz algebra $(\mathfrak{g},[\cdot,\cdot,\cdot]_{\mathfrak{g}})$, there is a natural adjoint representation   on  itself.
The corresponding maps $\rho^l, \rho^m$ and $\rho^r$ are given by
$$\rho^l(x,y,z)=\rho^m(x,y,z)=\rho^r(x,y,z)=[x,y,z]_{\mathfrak{g}}  ~\text{for all}~  x,y,z\in \mathfrak{g}.$$
\end{exam}

\begin{prop} \label{prop:3-Leibniz} \cite{XU}
Let $(\mathfrak{g},[\cdot,\cdot,\cdot]_{\mathfrak{g}})$ be a 3-Leibniz algebra, $V$ be a vector space, and let
$\rho^l:\mathfrak{g}\otimes \mathfrak{g}\otimes V\rightarrow V,
\rho^m:\mathfrak{g}\otimes V\otimes \mathfrak{g}\rightarrow V $ and $\rho^r:V\otimes \mathfrak{g}\otimes \mathfrak{g}\rightarrow V$ be trilinear maps.
Then $(V;\rho^l,\rho^m,\rho^r)$  will be   a representation of  $\mathfrak{g}$ if and only if $(\mathfrak{g}\oplus V, [\cdot,\cdot,\cdot]_{\ltimes})$ is a  3-Leibniz algebra,
where $[\cdot,\cdot,\cdot]_{\ltimes}$ is defined as
\begin{align*}
[(x,u),(y,v),(z,w)]_{\ltimes}=&([x,y,z]_{\mathfrak{g}},\rho^l(x,y,w)+\rho^m(x,v,z)+\rho^r(u,y,z)),
\end{align*}
for all $(x,u),(y,v),(z,w)\in \mathfrak{g}\oplus V$.
\end{prop}

In \cite{Loday} Loday introduced the notion of   dialgebras   in the study of
Leibniz algebras.  Inspired by his work, we  introduce the notion of  3-tri-Leibniz algebras.

\begin{definition}
A 3-tri-Leibniz algebra is a vector space $\mathfrak{g}$ endowed with   three  trilinear maps $[\cdot,\cdot,\cdot]_{\dashv},[\cdot,\cdot,\cdot]_{\perp},[\cdot,\cdot,\cdot]_{\vdash}:\mathfrak{g}\times\mathfrak{g}\times\mathfrak{g}\rightarrow\mathfrak{g}$, which adhere to the following equations:
\begin{align}
 &[a, b, [x, y, z]_{\triangle}]_{\dashv}=[[a, b, x]_{\dashv},y,z]_{\dashv}+ [x,  [a, b, y]_{\dashv},z]_{\perp}+ [x,y,[a, b, z]_{\dashv}]_{\vdash},\label{2.7}\\
 &[a, b, [x, y, z]_{\dashv}]_{\vdash}=[[a, b, x]_{\vdash},y,z]_{\dashv}+ [x,  [a, b, y]_{\triangle},z]_{\dashv}+ [x,y,[a, b, z]_{\triangle}]_{\dashv},\label{2.8}\\
  &[a, b, [x, y, z]_{\vdash}]_{\vdash}=[[a, b, x]_{\triangle},y,z]_{\vdash}+ [x,  [a, b, y]_{\triangle},z]_{\vdash}+ [x,y,[a, b, z]_{\vdash}]_{\vdash},\label{2.9}\\
    &[a, b, [x, y, z]_{\perp}]_{\vdash}=[[a, b, x]_{\triangle},y,z]_{\perp}+ [x,  [a, b, y]_{\vdash},z]_{\perp}+ [x,y,[a, b, z]_{\triangle}]_{\perp},\label{2.10}\\
 &[a, b, [x, y, z]_{\triangle}]_{\perp}=[[a, b, x]_{\perp},y,z]_{\dashv}+ [x,  [a, b, y]_{\perp},z]_{\perp}+ [x,y,[a, b, z]_{\perp}]_{\vdash},\label{2.11}
\end{align}
for all  $ x, y, z, a, b\in \mathfrak{g}$ and $\triangle\in \{\vdash,\dashv,\perp\}$.
\end{definition}

It is   important to note that the Eq. \eqref{2.9} says that the bracket $[\cdot,\cdot,\cdot]_{\vdash}$ is a 3-Leibniz bracket on
the vector space $\mathfrak{g}$.  However, the brackets $[\cdot,\cdot,\cdot]_{\dashv}$  and  $[\cdot,\cdot,\cdot]_{\perp}$  are not 3-Leibniz brackets on vector space $\mathfrak{g}$.

\begin{exam}
Any 3-Leibniz algebra $(\mathfrak{g},[\cdot, \cdot, \cdot]_\mathfrak{g})$ can be regarded as a 3-tri-Leibniz algebra in which $[\cdot,\cdot,\cdot]_{\dashv}=[\cdot,\cdot,\cdot]_{\perp}=[\cdot,\cdot,\cdot]_{\vdash}=[\cdot,\cdot, \cdot]_\mathfrak{g}$.
On the contrary, a 3-tri-Leibniz algebra, equipped with the same three operations, is merely a 3-Leibniz algebra.
\end{exam}

\begin{exam}
A differential 3-Leibniz algebra is a 3-Leibniz algebra $(\mathfrak{g},[\cdot, \cdot, \cdot]_\mathfrak{g})$  equipped with a linear map
$d: \mathfrak{g}\rightarrow \mathfrak{g}$ satisfying
$$d^2=0~~\text{and}~~ d[x, y, z]_\mathfrak{g}=[d(x), y, z]_\mathfrak{g}+[x, d(y), z]_\mathfrak{g}+[x, y, d(z)]_\mathfrak{g}  ~~\text{for all}~ x,y,z\in \mathfrak{g}.$$
Then it is easy to verify that $(\mathfrak{g},[\cdot,\cdot,\cdot]_{\dashv},[\cdot,\cdot,\cdot]_{\perp},[\cdot,\cdot,\cdot]_{\vdash})$ is a 3-tri-Leibniz algebra, where $[x,y,z]_{\dashv}=[x,y,d(z)]_\mathfrak{g},[x,y,z]_{\perp}=[x,d(y),z]_\mathfrak{g},[-,-,-]_{\vdash}=[d(x),y,z]_\mathfrak{g}$ for all $x,y,z\in\mathfrak{g}$.
\end{exam}

\begin{exam}
Let $(\mathfrak{g},[\cdot, \cdot, \cdot]_\mathfrak{g})$ be a 3-Leibniz algebra and  $(V;\rho^l,\rho^m,\rho^r)$  be a representation of  it.
Suppose $f:V\rightarrow \mathfrak{g}$ is a morphism between $\mathfrak{g}$-representations (from $V$ to the adjoint representation $\mathfrak{g}$), i.e.
 $$f(\rho^l(x,y,u))=[x, y, f(u)]_\mathfrak{g}, f(\rho^m(x,u,y))=[x, f(u), y]_\mathfrak{g}~~\text{and}~f(\rho^r(u,x,y))=[f(u), x, y]_\mathfrak{g}$$
for all~$x,y\in\mathfrak{g},u\in V. $
 Then there is a 3-tri-Leibniz algebra structure
on $V$ with the operations $[u,v,w]_{\dashv}=\rho^r(u,f(v),f(w)),[u,v,w]_{\perp}=\rho^m(f(u),v,f(w)),[u,v,w]_{\vdash}=\rho^l(f(u),f(v),w)$ for all $u,v,w\in V$.
\end{exam}

\begin{exam}
Let $(\mathfrak{g},[\cdot, \cdot, \cdot]_\mathfrak{g})$ be a 3-Leibniz algebra. Then the direct sum
$\underbrace{\mathfrak{g}\oplus\cdots\oplus\mathfrak{g}}_{n}$
carries a 3-tri-Leibniz algebra structure with the operations
\begin{align*}
[(x_1,\ldots,x_n),(y_1,\ldots, y_n),(z_1,\ldots,z_n)]_{\vdash}=&([x_1+\dots+x_n,y_1,z_1]_\mathfrak{g},\ldots,[x_1+\dots+x_n,y_n,z_n]_\mathfrak{g}),\\
[(x_1,\ldots,x_n),(y_1,\ldots, y_n),(z_1,\ldots,z_n)]_{\perp}=&([x_1, y_1+\dots+ y_n,z_1]_\mathfrak{g},\ldots,[x_n, y_1+\dots+ y_n,z_n]_\mathfrak{g}),\\
[(x_1,\ldots,x_n),(y_1,\ldots, y_n),(z_1,\ldots,z_n)]_{\dashv}=&([x_1, y_1, z_1+\dots+z_n]_\mathfrak{g},\ldots,[x_n, y_n, z_1+\dots+z_n]_\mathfrak{g}).
\end{align*}
\end{exam}

Another example  emerges from the representations of a 3-Leibniz algebra.

\begin{prop}\label{prop:3-tri-Leibniz algebra}
Let $(V;\rho^l,\rho^m,\rho^r)$  be a representation of a   3-Leibniz algebra $(\mathfrak{g},[\cdot,\cdot,\cdot]_{\mathfrak{g}})$.
Then $(\mathfrak{g}\oplus V, [\cdot,\cdot,\cdot]_{\dashv},[\cdot,\cdot,\cdot]_{\perp},[\cdot,\cdot,\cdot]_{\vdash})$ is a  3-tri-Leibniz algebra, where
\begin{align*}
[(x,u),(y,v),(z,w)]_{\vdash}=&([x,y,z]_{\mathfrak{g}},\rho^l(x,y,w)),\\
[(x,u),(y,v),(z,w)]_{\perp}=&([x,y,z]_{\mathfrak{g}},\rho^m(x,v,z)),\\
[(x,u),(y,v),(z,w)]_{\dashv}=&([x,y,z]_{\mathfrak{g}},\rho^r(u,y,z)),
\end{align*}
for all $(x,u),(y,v),(z,w)\in \mathfrak{g}\oplus V$.
This   3-tri-Leibniz algebra is called the hemisemidirect product  3-tri-Leibniz algebra and denoted by  $ \mathfrak{g}\ltimes_{l,m,r} V$.
\end{prop}

\begin{proof}
For any $(x,u),(y,v),(z,w),(a,s),(b,t)\in  \mathfrak{g}\oplus V$ and $\triangle\in \{\vdash,\dashv,\perp\}$, by Eqs.  \eqref{2.1}-\eqref{2.6}, we have
\begin{align*}
&[(a,s), (b,t), [(x,u), (y,v), (z,w)]_{\triangle}]_{\dashv}=([a,b,[x,y,z]_{\mathfrak{g}}]_{\mathfrak{g}},\rho^r(s, b,[x,y,z]_{\mathfrak{g}}))\\
&=\big([[a, b, x]_\mathfrak{g},y,z]_\mathfrak{g}+ [x,  [a, b, y]_\mathfrak{g},z]_\mathfrak{g}+ [x,y,[a, b, z]_\mathfrak{g}]_\mathfrak{g},\rho^r(\rho^r(s,b,x), y,z)+\\
&~~~~~ \rho^m(x,\rho^r(s,b,y),z)+\rho^l(x,y,\rho^r(s,b,z))\big)\\
&=\big([[a, b, x]_\mathfrak{g},y,z]_\mathfrak{g},\rho^r(\rho^r(s,b,x), y,z)\big)+\big([x,  [a, b, y]_\mathfrak{g},z]_\mathfrak{g},\rho^m(x,\rho^r(s,b,y),z)\big)+\\
&~~~~~\big([x,y,[a, b, z]_\mathfrak{g}]_\mathfrak{g},\rho^l(x,y,\rho^r(s,b,z))\big)\\
&=[[(a,s), (b,t), (x,u)]_{\dashv}, (y,v), (z,w)]_{\dashv}+[(x,u),[(a,s), (b,t),  (y,v)]_{\dashv}, (z,w)]_{\perp}+\\
&~~~~~[(x,u), (y,v), [(a,s), (b,t),  (z,w)]_{\dashv}]_{\vdash}.
\end{align*}
Similarly, we get
\begin{align*}
&[(a,s), (b,t), [(x,u), (y,v), (z,w)]_{\dashv}]_{\vdash}=([a,b,[x,y,z]_{\mathfrak{g}}]_{\mathfrak{g}},\rho^l(a, b,\rho^r(u,y,z)))\\
&=\big([[a, b, x]_\mathfrak{g},y,z]_\mathfrak{g}+ [x,  [a, b, y]_\mathfrak{g},z]_\mathfrak{g}+ [x,y,[a, b, z]_\mathfrak{g}]_\mathfrak{g},\rho^r(\rho^l(a,b,u), y,z)+ \\
&~~~~~\rho^r(u,[a,b,y]_{\mathfrak{g}},z)+\rho^r(u,y,[a,b,z]_{\mathfrak{g}})\big)\\
&=[[(a,s), (b,t), (x,u)]_{\vdash}, (y,v), (z,w)]_{\dashv}+[(x,u),[(a,s), (b,t),  (y,v)]_{\triangle}, (z,w)]_{\dashv}+\\
&~~~~~[(x,u), (y,v), [(a,s), (b,t),  (z,w)]_{\triangle}]_{\dashv}.
\end{align*}
and
\begin{align*}
&[(a,s), (b,t), [(x,u), (y,v), (z,w)]_{\vdash}]_{\vdash}=([a,b,[x,y,z]_{\mathfrak{g}}]_{\mathfrak{g}},\rho^l(a, b,\rho^l(x,y,w)))\\
&=\big([[a, b, x]_\mathfrak{g},y,z]_\mathfrak{g}+ [x,  [a, b, y]_\mathfrak{g},z]_\mathfrak{g}+ [x,y,[a, b, z]_\mathfrak{g}]_\mathfrak{g},\rho^l([a,b,x]_{\mathfrak{g}}, y, w)+\\
&~~~~~ \rho^l(x,[a,b,y]_{\mathfrak{g}},w)+\rho^l(x,y,\rho^l(a,b,w))\big)\\
&=[[(a,s), (b,t), (x,u)]_{\triangle}, (y,v), (z,w)]_{\vdash}+[(x,u),[(a,s), (b,t),  (y,v)]_{\triangle}, (z,w)]_{\vdash}+\\
&~~~~~[(x,u), (y,v), [(a,s), (b,t),  (z,w)]_{\vdash}]_{\vdash}.
\end{align*}
Thus, the  Eqs.  \eqref{2.7}-\eqref{2.9} hold. Finally,
\begin{align*}
&[(a,s), (b,t), [(x,u), (y,v), (z,w)]_{\perp}]_{\vdash}=([a,b,[x,y,z]_{\mathfrak{g}}]_{\mathfrak{g}},\rho^l(a, b,\rho^m(x,v,z)))\\
&=\big([[a, b, x]_\mathfrak{g},y,z]_\mathfrak{g}+ [x,  [a, b, y]_\mathfrak{g},z]_\mathfrak{g}+ [x,y,[a, b, z]_\mathfrak{g}]_\mathfrak{g},\rho^m([a,b,x]_{\mathfrak{g}}, v, z)+\\
&~~~~~\rho^m(x,\rho^l(a,b,v),z)+\rho^m(x,v,[a,b,z]_{\mathfrak{g}})\big)\\
&=[[(a,s), (b,t), (x,u)]_{\triangle}, (y,v), (z,w)]_{\perp}+[(x,u),[(a,s), (b,t),  (y,v)]_{\vdash}, (z,w)]_{\perp}+\\
&~~~~~[(x,u), (y,v), [(a,s), (b,t),  (z,w)]_{\triangle}]_{\perp}.
\end{align*}
and
\begin{align*}
&[(a,s), (b,t), [(x,u), (y,v), (z,w)]_{\triangle}]_{\perp}=([a,b,[x,y,z]_{\mathfrak{g}}]_{\mathfrak{g}},\rho^m(a, t,[x,y,z]_{\mathfrak{g}}))\\
&=\big([[a, b, x]_\mathfrak{g},y,z]_\mathfrak{g}+ [x,  [a, b, y]_\mathfrak{g},z]_\mathfrak{g}+ [x,y,[a, b, z]_\mathfrak{g}]_\mathfrak{g},\rho^r(\rho^m(a,t,x), y,z)+\\
&~~~~~ \rho^m(x,\rho^m(a,t,y),z)+\rho^l(x,y,\rho^m(a,t,z))\big)\\
&=[[(a,s), (b,t), (x,u)]_{\perp}, (y,v), (z,w)]_{\dashv}+[(x,u),[(a,s), (b,t),  (y,v)]_{\perp}, (z,w)]_{\perp}+\\
&~~~~~[(x,u), (y,v), [(a,s), (b,t),  (z,w)]_{\perp}]_{\vdash}.
\end{align*}
So the  Eqs. \eqref{2.10} and \eqref{2.11} hold  and we complete the proof.
\end{proof}

\section{Embedding tensors on 3-Leibniz algebras }\label{sec: Embedding tensors}
\def\theequation{\arabic{section}.\arabic{equation}}
\setcounter{equation} {0}

In this section, we first introduce embedding tensors  and  averaging operators on 3-Leibniz algebras.
Then we give two equivalent characterization theorems via graphs and Nijenhuis operators.
Following this, we show that an embedding tensor
 or an   averaging operator induces a  3-tri-Leibniz algebra. Finally, we show that any
3-tri-Leibniz algebra can be embedded into an averaging 3-Leibniz algebra.

\begin{definition}
Let $(\mathfrak{g},[\cdot, \cdot, \cdot]_\mathfrak{g})$ be a 3-Leibniz algebra and  $(V;\rho^l,\rho^m,\rho^r)$  be a representation of  it.  An embedding tensor on $\mathfrak{g}$ with respect to the representation $V$
is a linear map
$T:V\rightarrow  \mathfrak{g}$ that satisfies
\begin{align}
 [Tu,Tv,Tw]_\mathfrak{g}=T\rho^l(Tu,Tv,w)=T\rho^m(Tu,v,Tw)=T\rho^r(u,Tv,Tw) \label{3.1}
\end{align}
for all  $u,v,w\in V$.
\end{definition}

The embedding tensor of a 3-Leibniz algebra $(\mathfrak{g},[\cdot, \cdot, \cdot]_\mathfrak{g})$  with respect to the adjoint representation
$(\mathfrak{g};\rho^l,\rho^m,\rho^r)$ is called an averaging operator. In this case the Eq. \eqref{3.1} can be written as
\begin{align}
 [Tx,Ty,Tz]_\mathfrak{g}= T[Tx,Ty,z]_\mathfrak{g}= T[Tx,y,Tz]_\mathfrak{g}=T [x,Ty,Tz]_\mathfrak{g} \label{3.2}
\end{align}
for all  $x,y,z\in  \mathfrak{g}$.

\begin{definition}
An averaging 3-Leibniz algebra is a pair $(\mathfrak{g},T)$ consisting of a 3-Leibniz algebra $\mathfrak{g}$ endowed with an
averaging operator $T$ on it.
\end{definition}

\begin{exam}
Let  $(\mathfrak{g},[\cdot, \cdot, \cdot]_\mathfrak{g})$ be a 3-Leibniz algebra . Then the identity map $id: \mathfrak{g}\rightarrow \mathfrak{g}$ is an averaging
operator on $\mathfrak{g}$.
\end{exam}

\begin{exam}
Let $(\mathfrak{g},[\cdot, \cdot,\cdot]_\mathfrak{g},d)$ be a differential 3-Leibniz algebra.
Then it is easy to see that  $d: \mathfrak{g}\rightarrow \mathfrak{g}$ is an averaging
operator on $\mathfrak{g}$.
\end{exam}

\begin{exam}
Let $(\mathfrak{g},[\cdot, \cdot, \cdot]_\mathfrak{g})$ be a 3-Leibniz algebra and  $(V;\rho^l,\rho^m,\rho^r)$  be a representation of  it.
Suppose $f:V\rightarrow \mathfrak{g}$ is a morphism of $\mathfrak{g}$-representations  from $V$ to the adjoint representation $\mathfrak{g}$.
 Then $f$ is an   embedding tensor on $\mathfrak{g}$ with respect to the representation $V$.
\end{exam}

\begin{exam}
Let $(\mathfrak{g},[\cdot, \cdot, \cdot]_\mathfrak{g})$ be a 3-Leibniz algebra. Then the space $\underbrace{\mathfrak{g}\oplus\cdots\oplus\mathfrak{g}}_{n}$
can be given a representation
of the 3-Leibniz algebra $\mathfrak{g}$ with the three actions given by
\begin{align*}
\rho^l(x,y,(z_1,\ldots,z_n))=&([x,y,z_1]_\mathfrak{g},\ldots,[x,y,z_n]_\mathfrak{g}),\\
\rho^m(x,(z_1,\ldots, z_n),y)=&([x, z_1,y]_\mathfrak{g},\ldots,[x,  z_n,y]_\mathfrak{g}),\\
\rho^r((z_1,\ldots,z_n),x,y)=&([z_1, x,y]_\mathfrak{g},\ldots,[z_n, x,y]_\mathfrak{g})
\end{align*}
for all $x,y\in \mathfrak{g}$ and $(z_1,\ldots,z_n)\in \mathfrak{g}\oplus\cdots\oplus\mathfrak{g}.$
With this notation, the map
$$T:\mathfrak{g}\oplus\cdots\oplus\mathfrak{g}\rightarrow\mathfrak{g},(z_1,\ldots,z_n)\mapsto z_1+\dots+z_n$$
is an   embedding tensor on $\mathfrak{g}$ with respect to the representation $\mathfrak{g}\oplus\cdots\oplus\mathfrak{g}$. Moreover, for each $1\leq i\leq n$, the $i$-th projection map $$T_i:\mathfrak{g}\oplus\cdots\oplus\mathfrak{g}\rightarrow\mathfrak{g},(z_1,\ldots,z_n)\mapsto z_i$$
is also an   embedding tensor on $\mathfrak{g}$ with respect to the representation $\mathfrak{g}\oplus\cdots\oplus\mathfrak{g}$.
\end{exam}

\begin{definition}
Let  $(\mathfrak{h},[\cdot, \cdot, \cdot]_\vdash, [\cdot, \cdot, \cdot]_\dashv, [\cdot, \cdot, \cdot]_\perp)$   be a 3-tri-Leibniz algebra.
A Nejinhuis operator on   $\mathfrak{h}$ is a linear map $N:\mathfrak{h}\rightarrow \mathfrak{h}$, such that
 \begin{align*}
[Nx,Ny,Nz]_\triangle=&N\big([x,Ny,Nz]_\triangle+[Nx,y,Nz]_\triangle+[Nx,Ny,z]_\triangle)-\\
&N^2([Nx,y,z]_\triangle+[x,Ny,z]_\triangle+[x,y,Nz]_\triangle)+N^3[x,y,z]_\triangle
\end{align*}
for all $x,y,z\in \mathfrak{h}$ and $\triangle\in\{\vdash,\dashv,\perp\}$.
\end{definition}

Next, we present two equivalent characterization of an embedding tensor in terms of Nijenhuis operators and  its associated graph.

\begin{theorem} \label{theorem:Nijenhuis}
A linear map $T:V\rightarrow  \mathfrak{g}$ is an embedding tensor   on a 3-Leibniz algebra $(\mathfrak{g},[\cdot, \cdot, \cdot]_\mathfrak{g})$ with respect to the
representation  $(V;\rho^l,\rho^m,\rho^r)$  if and only if the map
 $$N_T:\mathfrak{g}\oplus V\rightarrow \mathfrak{g}\oplus V,(x,u)\mapsto(Tu,0)$$
is a Nijenhuis operator on the   hemisemidirect product
3-tri-Leibniz algebra $ \mathfrak{g}\ltimes_{l,m,r} V$.
\end{theorem}

\begin{proof}
Let $T:V\rightarrow  \mathfrak{g}$  be a linear map.
For any $(x,u),(y,v),(z,w)\in \mathfrak{g}\oplus V$,  we have
\begin{align*}
&[N_T(x,u),N_T(y,v),N_T(z,w)]_{\vdash}-N_T\big([(x,u),N_T(y,v),N_T(z,w)]_{\vdash}+[N_T(x,u),(y,v),N_T(z,w)]_{\vdash}+\\
&[N_T(x,u),N_T(y,v),(z,w)]_{\vdash}\big)+N_T^2\big([N_T(x,u),(y,v),(z,w)]_{\vdash}+[(x,u),N_T(y,v),(z,w)]_{\vdash}+\\
&[(x,u),(y,v),N_T(z,w)]_{\vdash}\big)-N_T^3[(x,u),(y,v),(z,w)]_{\vdash}\\
&=([Tu,Tv,Tw]_\mathfrak{g},0)-N_T\big(([x,Tv,Tw]_\mathfrak{g},0)+([Tu,y,Tw]_\mathfrak{g},0)+([Tu,Tv,z]_\mathfrak{g},\rho^l(Tu,Tv,w))\big)+\\
&~~~~N_T^2\big(([Tu,y,z]_\mathfrak{g},\rho^l(Tu,y,w))+([x,Tv,z]_\mathfrak{g},\rho^l(x,Tv,w))+([x,y,Tw]_\mathfrak{g},0)\big)+\\
&~~~~N_T^3([x,y,z]_\mathfrak{g},\rho^l(x,y,w))\\
&=([Tu,Tv,Tw]_\mathfrak{g},0)-(T\rho^l(Tu,Tv,w),0).
\end{align*}
Similarly, we obtain
\begin{align*}
&[N_T(x,u),N_T(y,v),N_T(z,w)]_{\dashv}-N_T\big([(x,u),N_T(y,v),N_T(z,w)]_{\dashv}+[N_T(x,u),(y,v),N_T(z,w)]_{\dashv}+\\
&[N_T(x,u),N_T(y,v),(z,w)]_{\dashv}\big)+N_T^2\big([N_T(x,u),(y,v),(z,w)]_{\dashv}+[(x,u),N_T(y,v),(z,w)]_{\dashv}+\\
&[(x,u),(y,v),N_T(z,w)]_{\dashv}\big)-N_T^3[(x,u),(y,v),(z,w)]_{\dashv}\\
&=([Tu,Tv,Tw]_\mathfrak{g},0)-(T\rho^r(u,Tv,Tw),0).
\end{align*}
and
\begin{align*}
&[N_T(x,u),N_T(y,v),N_T(z,w)]_{\perp}-N_T\big([(x,u),N_T(y,v),N_T(z,w)]_{\perp}+[N_T(x,u),(y,v),N_T(z,w)]_{\perp}+\\
&[N_T(x,u),N_T(y,v),(z,w)]_{\perp}\big)+N_T^2\big([N_T(x,u),(y,v),(z,w)]_{\perp}+[(x,u),N_T(y,v),(z,w)]_{\perp}+\\
&[(x,u),(y,v),N_T(z,w)]_{\perp}\big)-N_T^3[(x,u),(y,v),(z,w)]_{\perp}\\
&=([Tu,Tv,Tw]_\mathfrak{g},0)-(T\rho^m(Tu,v,Tw),0).
\end{align*}
Therefore $N_T$ is a Nijenhuis operator on the hemisemidirect product  3-tri-Leibniz algebra  $ \mathfrak{g}\ltimes_{l,m,r} V$ if and only if   $T$ is an embedding tensor.
\end{proof}

\begin{theorem} \label{theorem:graph}
A linear map $T:V\rightarrow  \mathfrak{g}$ is an embedding tensor   on a 3-Leibniz algebra $(\mathfrak{g},[\cdot, \cdot, \cdot]_\mathfrak{g})$ with respect to the
representation  $(V;\rho^l,\rho^m,\rho^r)$  if and only if the graph $Gr(T)=\{(Tu,u)~|~u\in V\}$ is a subalgebra of the hemisemidirect product
3-tri-Leibniz algebra $ \mathfrak{g}\ltimes_{l,m,r} V$.
\end{theorem}

\begin{proof}
Let $T:V\rightarrow  \mathfrak{g}$  be a linear map.
For any $u,v,w\in V$,  we have
\begin{align*}
[(Tu,u),(Tv,v),(Tw,w)]_{\vdash}=&([Tu,Tv,Tw]_{\mathfrak{g}},\rho^l(Tu,Tv,w)),\\
[(Tu,u),(Tv,v),(Tw,w)]_{\perp}=&([Tu,Tv,Tw]_{\mathfrak{g}},\rho^m(Tu,v,Tw)),\\
[(Tu,u),(Tv,v),(Tw,w)]_{\dashv}=&([Tu,Tv,Tw]_{\mathfrak{g}},\rho^r(u,Tv,Tw)),
\end{align*}
The above three elements are in $Gr(T)$ if and only if $[Tu,Tv,Tw]_{\mathfrak{g}}=T\rho^l(Tu,Tv,w)$, $[Tu,Tv,Tw]_{\mathfrak{g}}=T\rho^m(Tu,v,Tw)$ and $[Tu,Tv,Tw]_{\mathfrak{g}}=T\rho^r(u,Tv,Tw)$.
Therefore,  the graph $Gr(T)$ is a subalgebra of the hemisemidirect product  3-tri-Leibniz algebra  $ \mathfrak{g}\ltimes_{l,m,r} V$ if and only if   $T$ is an embedding tensor.
\end{proof}

Since  $Gr(T)$ is linearly isomorphic to $V$,  the aforementioned theorem leads to the following outcome.

\begin{prop} \label{prop:ET}
Let $T:V\rightarrow  \mathfrak{g}$ be an embedding tensor on a 3-Leibniz algebra $(\mathfrak{g},[\cdot, \cdot, \cdot]_\mathfrak{g})$ with respect to the
representation  $(V;\rho^l,\rho^m,\rho^r)$.  Then the representation space $V$ inherits a 3-tri-Leibniz algebra structure with the operations
\begin{align*}
[u,v,w]_{\vdash}^T=\rho^l(Tu,Tv,w),~~[u,v,w]_{\perp}^T=\rho^m(Tu,v,Tw)~\text{and}~[u,v,w]_{\dashv}^T=\rho^r(u,Tv,Tw)
\end{align*}
for  all $u,v,w\in V$.   Moreover, $T$ is a homomorphism from the 3-tri-Leibniz algebra $(V,[\cdot,\cdot,\cdot]_{\vdash}^T,[\cdot,\cdot,\cdot]_{\perp}^T,[\cdot,\cdot,\cdot]_{\dashv}^T)$ to the 3-Leibniz algebra $(\mathfrak{g},[\cdot,\cdot,\cdot]_\mathfrak{g})$.
\end{prop}

\begin{coro}
Let $(\mathfrak{g},[\cdot, \cdot,\cdot]_\mathfrak{g})$ be a 3-Leibniz algebra and $T: \mathfrak{g}\rightarrow  \mathfrak{g}$  be an averaging operator. Then,
$\mathfrak{g}$ is also equipped with  a 3-tri-Leibniz algebra structure defined by
 \begin{align*}
[x,y,z]_{\vdash}^T=[Tx,Ty,z]_\mathfrak{g},~~[x,y,z]_{\perp}^T=[Tx,y,Tz]_\mathfrak{g}~\text{and}~[x,y,z]_{\dashv}^T=[x,Ty,Tz]_\mathfrak{g}.
\end{align*}
\end{coro}

The above result  indicates that  an embedding tensor  on a 3-Leibniz algebra with respect to a representation induces a 3-tri-Leibniz algebra structure on the underlying representation space.
 The subsequent result establishes the converse relationship.

\begin{theorem}
Every 3-tri-Leibniz algebra is induced by an   embedding tensor on a 3-Leibniz algebra
with respect to a representation.
\end{theorem}

\begin{proof}
Let  $(\mathfrak{g},[\cdot, \cdot, \cdot]_\vdash, [\cdot, \cdot, \cdot]_\dashv, [\cdot, \cdot, \cdot]_\perp)$   be a 3-tri-Leibniz algebra.
The vector subspace $I_{\mathfrak{g}}$ spanned by
\begin{align*}
\big\{[x, y, z]_\vdash- [x, y, z]_\dashv,~[x, y, z]_\dashv- [x, y, z]_\perp,~[x, y, z]_\vdash-   [x, y, z]_\perp~:~x,y,z\in\mathfrak{g}\big\}
\end{align*}
is an ideal of $\mathfrak{g}$.
 Then, the quotient algebra $\mathfrak{g}_{\mathrm{3Leib}}=\frac{\mathfrak{g}}{I_{\mathfrak{g}}}$ carries a
3-Leibniz algebra structure with the bracket
$$[\overline{x},\overline{y},\overline{z}]=\overline{[x, y, z]_\vdash}=\overline{[x, y, z]_\dashv}= \overline{[x,y, z]_\perp}$$
for all $x,y,z\in \mathfrak{g}$, here $\overline{x}$ denotes the class of an element $x\in \mathfrak{g}$.
We define three trilinear maps
$\rho^l: \mathfrak{g}_{\mathrm{3Leib}}\times\mathfrak{g}_{\mathrm{3Leib}}\times\mathfrak{g}\rightarrow\mathfrak{g}$, $\rho^m: \mathfrak{g}_{\mathrm{3Leib}}\times\mathfrak{g}\times\mathfrak{g}_{\mathrm{3Leib}}\rightarrow\mathfrak{g}$ and
$\rho^r: \mathfrak{g}\times\mathfrak{g}_{\mathrm{3Leib}}\times\mathfrak{g}_{\mathrm{3Leib}}\rightarrow\mathfrak{g}$
by
\begin{align*}
\rho^l(\overline{x},\overline{y},z)=[x, y, z]_\vdash,~~\rho^m(\overline{x},z, \overline{y})=[x, z, y]_\perp,~~\rho^r(z,\overline{x},  \overline{y})=[z,x,   y]_\dashv
\end{align*}
for all $\overline{x},\overline{y}\in \mathfrak{g}_{\mathrm{3Leib}}$ and $z\in \mathfrak{g}$.
It is easy to verify that the maps $\rho^l,\rho^m,\rho^r$
define a representation of the 3-Leibniz algebra $\mathfrak{g}_{\mathrm{3Leib}}$ on the vector
space  $\mathfrak{g}$ . Moreover, the quotient map $T:  \mathfrak{g}\rightarrow  \mathfrak{g}_{\mathrm{3Leib}}, x\mapsto \overline{x}$ is an   embedding tensor as
\begin{align*}
[Tx,Ty,Tz]=[\overline{x},\overline{y},\overline{z}]=
\left\{ \begin{array}{lll}
=\overline{[x, y, z]_\vdash}=\overline{\rho^l(\overline{x},\overline{y},z)}=T\rho^l(Tx,Ty,z),\\
=\overline{[x, y, z]_\dashv}=\overline{\rho^r(x,\overline{y},  \overline{z})}=T\rho^r(x,Ty,Tz),\\
=\overline{[x,y, z]_\perp}=\overline{\rho^m(\overline{x},y,  \overline{z})}=T\rho^r(Tx,y,Tz)
 \end{array}  \right.
 \end{align*}
 for all $x,y,z\in \mathfrak{g}$. Let $(\mathfrak{g},[\cdot, \cdot, \cdot]_\vdash^T, [\cdot, \cdot, \cdot]_\dashv^T, [\cdot, \cdot, \cdot]_\perp^T)$  be the 3-tri-Leibniz algebra structure on $\mathfrak{g}$ induced by the  embedding tensor $T$.
  Then we have
\begin{align*}
[x, y, z]_\vdash^T=\rho^l(\overline{x},\overline{y},z), ~~[x,y,z]_{\perp}^T=\rho^m(\overline{x},y,\overline{z})~\text{and}~[x,y,z]_{\dashv}^T=\rho^r(x,\overline{y},\overline{z}).
\end{align*}
  This demonstrates that the induced 3-tri-Leibniz algebra $(\mathfrak{g},[\cdot,\cdot, \cdot]_\vdash^T, [\cdot, \cdot,\cdot]_\dashv^T, [\cdot, \cdot, \cdot]_\perp^T)$   is identical to the specified one.
\end{proof}

\begin{remark}
It is crucial to note that an arbitrary 3-tri-Leibniz algebra may not be induced from an averaging 3-Leibniz algebra.
However, any 3-tri-Leibniz algebra can be embedded into an averaging 3-Leibniz
algebra.
More precisely, let   $(\mathfrak{g},[\cdot, \cdot, \cdot]_\vdash, [\cdot, \cdot, \cdot]_\dashv, $ $[\cdot, \cdot, \cdot]_\perp)$   be a 3-tri-Leibniz algebra.  Then by Proposition \ref{prop:3-Leibniz},
the direct sum $\mathfrak{g}_{\mathrm{3Leib}}\oplus \mathfrak{g}$
inherits a 3-Leibniz algebra structure with the bracket
\begin{align*}
[(\overline{x},x),(\overline{y},y),(\overline{z},z)]_{\ltimes}=([\overline{x},\overline{y},\overline{z}],\rho^l(\overline{x},\overline{y},z)+\rho^m(\overline{x},y,\overline{z})+\rho^r(x,\overline{y},\overline{z}))
 \end{align*}
 for all $(\overline{x},x),(\overline{y},y),(\overline{z},z)\in \mathfrak{g}_{\mathrm{3Leib}}\oplus \mathfrak{g}.$
 Moreover, the map
 $$T:\mathfrak{g}_{\mathrm{3Leib}}\oplus \mathfrak{g}\rightarrow\mathfrak{g}_{\mathrm{3Leib}}\oplus \mathfrak{g},~~(\overline{x},y)\mapsto (\overline{y},0)$$
  is an averaging operator. Thus, $(\mathfrak{g}_{\mathrm{3Leib}}\oplus \mathfrak{g},T)$ is an averaging 3-Leibniz algebra.
  Then the inclusion map
  $$i:\mathfrak{g}\rightarrow\mathfrak{g}_{\mathrm{3Leib}}\oplus \mathfrak{g},~~ x\mapsto(0,x)$$
  is an embedding of the 3-tri-Leibniz algebra $(\mathfrak{g},[\cdot, \cdot, \cdot]_\vdash, [\cdot, \cdot,\cdot]_\dashv, [\cdot, \cdot, \cdot]_\perp)$ into the
averaging 3-Leibniz algebra $(\mathfrak{g}_{\mathrm{3Leib}}\oplus \mathfrak{g},T)$.
\end{remark}

\section{ 3-tri-Leibniz dialgebras }\label{sec: 3-tri-Leibniz dialgebras}
\def\theequation{\arabic{section}.\arabic{equation}}
\setcounter{equation} {0}
This section duplicates a 3-tri-Leibniz algebra, resulting in what we will refer to as a 3-tri-Leibniz dialgebra.
For this purpose, we  first introduce homomorphic embedding tensors on a 3-Leibniz algebra with respect to a   action.

\begin{definition}
Let  $(\mathfrak{g},[\cdot, \cdot,\cdot]_\mathfrak{g})$ and $(\mathfrak{h},[\cdot, \cdot,\cdot]_\mathfrak{h})$ be two 3-Leibniz algebras. We say that $\mathfrak{g}$ acts on $\mathfrak{h}$
if there  are three trilinear maps:
\begin{align*}
\rho^l:\mathfrak{g}\otimes \mathfrak{g}\otimes \mathfrak{h}\rightarrow \mathfrak{h},~~
\rho^m:\mathfrak{g}\otimes \mathfrak{h}\otimes \mathfrak{g}\rightarrow \mathfrak{h}, ~~
\rho^r:\mathfrak{h}\otimes \mathfrak{g}\otimes \mathfrak{g}\rightarrow \mathfrak{h},
\end{align*}
such that $(\mathfrak{h};\rho^l,\rho^m,\rho^r)$ is a representation of $\mathfrak{g}$ and for any
$a,b,x,y,z\in\mathfrak{g}$ and $s,u,v,w\in \mathfrak{h}$,
\begin{align}
&\rho^l\big(a,b,[u,v,w]_\mathfrak{h}\big)=[\rho^l(a,b,u),v,w]_\mathfrak{h}+[u, \rho^l(a,b,v),w]_\mathfrak{h}+[u,v,\rho^l(a,b,w)]_\mathfrak{h},\label{4.1}\\
&[\rho^m(a,s,x),v,w]_\mathfrak{h}+[u, \rho^m(a,s,y),w]_\mathfrak{h}+[u,v,\rho^m(a,s,z)]_\mathfrak{h}=0,\label{4.2}\\
&[\rho^r(s,a,x),v,w]_\mathfrak{h}+[u, \rho^r(s,a,y),w]_\mathfrak{h}+[u,v,\rho^r(s,a,z)]_\mathfrak{h}=0.\label{4.3}\\
&\rho^l\big(a,b,[u,v,w]_\mathfrak{h}\big)=[u,v,\rho^l(a,b,w)]_{\mathfrak{h}},\label{4.4}\\
&\rho^m(a,[u,v,w]_\mathfrak{h},b)=[u,v,\rho^m(a,w,b)]_{\mathfrak{h}},\label{4.5}\\
&\rho^r([u,v,w]_\mathfrak{h},a,b)=[u,v,\rho^r(w,a,b)]_{\mathfrak{h}}.\label{4.6}
\end{align}
\end{definition}

\begin{prop}
Let  $(\mathfrak{g},[\cdot, \cdot,\cdot]_\mathfrak{g})$ and $(\mathfrak{h},[\cdot, \cdot,\cdot]_\mathfrak{h})$ be two 3-Leibniz algebras,  and let
$\rho^l:\mathfrak{g}\otimes \mathfrak{g}\otimes \mathfrak{h}\rightarrow \mathfrak{h},
\rho^m:\mathfrak{g}\otimes \mathfrak{h}\otimes \mathfrak{g}\rightarrow \mathfrak{h} $ and $\rho^r:\mathfrak{h}\otimes \mathfrak{g}\otimes \mathfrak{g}\rightarrow \mathfrak{h}$ be trilinear maps.
Then, the tuple $(\mathfrak{h},[\cdot, \cdot, \cdot]_\mathfrak{h};\rho^l,\rho^m,\rho^r)$ is an action of $\mathfrak{g}$ if and only if
$\mathfrak{g}\oplus \mathfrak{h}$ carries a 3-Leibniz algebra structure with operation given by
\begin{align*}
 [(x,u),(y,v),(z,w)]_{\bowtie}=([x,y,z]_\mathfrak{g},\rho^l(x,y,w)+\rho^m(x,v,z)+\rho^r(u,y,z)+[u,v,w]_\mathfrak{h})
\end{align*}
for all $(x,u),(y,v),(z,w)\in \mathfrak{g}\oplus\mathfrak{h}$, which is called the semi-direct product of $\mathfrak{g}$ with $\mathfrak{h}$.
\end{prop}

\begin{proof}
For any $(a,s), (b,t), (x,u),(y,v),(z,w)\in \mathfrak{g}\oplus\mathfrak{h}$, we have
\begin{align*}
&[(a,s), (b,t),  [(x,u),(y,v),(z,w)]_{\bowtie}]_{\bowtie}-[[(a,s), (b,t),  (x,u)]_{\bowtie},(y,v),(z,w)]_{\bowtie}-\\
&[(x,u),[(a,s), (b,t),  (y,v)]_{\bowtie},(z,w)]_{\bowtie}-[(x,u),(y,v), [(a,s), (b,t),  (z,w)]_{\bowtie}]_{\bowtie}\\
&=\Big([a, b,[x,y,z]_\mathfrak{g}]_\mathfrak{g}-[[a,b,x]_\mathfrak{g},y,z]_{ \mathfrak{g}}-[x,[a,b,y]_\mathfrak{g},z]_{\mathfrak{g}}-[x,y,[a,b,z]_\mathfrak{g}]_\mathfrak{g},\\
&\rho^l\big(a,b,\rho^l(x,y,w)+\rho^m(x,v,z)+\rho^r(u,y,z)+[u,v,w]_\mathfrak{h}\big)+\rho^m(a,t,[x,y,z]_\mathfrak{g})+\\
&\rho^r(s,b,[x,y,z]_\mathfrak{g})+[s,t,\rho^l(x,y,w)+\rho^m(x,v,z)+\rho^r(u,y,z)+[u,v,w]_\mathfrak{h}]_{\mathfrak{h}}-\\
&\rho^l([a,b,x]_\mathfrak{g},y,w)-\rho^m([a,b,x]_\mathfrak{g},v,z)-\rho^r(\rho^l(a,b,u)+\rho^m(a,t,x)+\rho^r(s,b,x)+[s,t,u]_\mathfrak{h},y,z)-\\
&[\rho^l(a,b,u)+\rho^m(a,t,x)+\rho^r(s,b,x)+[s,t,u]_\mathfrak{h},v,w]_\mathfrak{h}-\\
&\rho^l(x,[a,b,y]_\mathfrak{g},w)-\rho^m(x,\rho^l(a,b,v)+\rho^m(a,t,y)+\rho^r(s,b,y)+[s,t,v]_\mathfrak{h},z)-\\
&\rho^r(u,[a,b,y]_\mathfrak{g},z)-[u, \rho^l(a,b,v)+\rho^m(a,t,y)+\rho^r(s,b,y)+[s,t,v]_\mathfrak{h},w]_\mathfrak{h}-\\
&\rho^l\big(x,y,\rho^l(a,b,w)+\rho^m(a,t,z)+\rho^r(s,b,z)+[s,t,w]_\mathfrak{h}\big)-\rho^m(x,v,[a,b,z]_\mathfrak{g})-\\
&\rho^r(u,y,[a,b,z]_\mathfrak{g})-[u,v,\rho^l(a,b,w)+\rho^m(a,t,z)+\rho^r(s,b,z)+[s,t,w]_\mathfrak{h}]_\mathfrak{h}\Big).
\end{align*}
Then, $(\mathfrak{g}\oplus \mathfrak{h},  [\cdot,\cdot,\cdot]_{\bowtie})$ is a 3-Leibniz algebra if and only if  $(\mathfrak{h};\rho^l,\rho^m,\rho^r)$ is a representation of $\mathfrak{g}$ and
Eqs. \eqref{4.1}-\eqref{4.6} are satisfied.
That is $(\mathfrak{g}\oplus \mathfrak{h},  [\cdot,\cdot,\cdot]_{\bowtie})$ is a 3-Leibniz algebra if and only if  $(\mathfrak{h},[\cdot, \cdot, \cdot]_\mathfrak{h};\rho^l,\rho^m,\rho^r)$ is an action of   $(\mathfrak{g},[\cdot, \cdot,\cdot]_\mathfrak{g})$.
\end{proof}

In the following, we introduce the notion of homomorphic  embedding tensor on a 3-Leibniz algebra with respect to a given action.

\begin{definition}
 A linear map $T:\mathfrak{h}\rightarrow \mathfrak{g}$ is called a homomorphic embedding tensor on the
 3-Leibniz algebra  $(\mathfrak{g},[\cdot, \cdot, \cdot]_\mathfrak{g})$  with respect to a   action $(\mathfrak{h},[\cdot, \cdot, \cdot]_\mathfrak{h};\rho^l,\rho^m,\rho^r)$  if $T$ qualifies both as an embedding tensor and a 3-Leibniz algebra homomorphism.
\end{definition}

A homomorphic embedding tensor on a
 3-Leibniz algebra  $\mathfrak{g}$  with respect to the adjoint representation $(\mathfrak{g}; \rho^l,\rho^m,\rho^r)$
is called a homomorphic averaging operator.

\begin{exam}
A crossed module of  3-Leibniz algebras  is a tuple $(\mathfrak{g},\mathfrak{h}, T,  \rho^l,\rho^m,\rho^r)$ in which $(\mathfrak{g},[\cdot, \cdot, \cdot]_\mathfrak{g})$
is a  3-Leibniz algebra, $(\mathfrak{h},[\cdot, \cdot,\cdot]_\mathfrak{h};\rho^l,\rho^m,\rho^r)$ is an action of   $\mathfrak{g}$ and $T:\mathfrak{h}\rightarrow \mathfrak{g}$  is a 3-Leibniz algebra homomorphism,
where
\begin{align*}
&T\rho^l(x,y,u)=[x,y,Tu]_\mathfrak{g},~~ T\rho^m(x,u,y)=[x,Tu,y]_\mathfrak{g},~~T\rho^r(u,x,y)=[Tu,x,y]_\mathfrak{g},\\
&\rho^l(Tu,Tv,w)=\rho^m(Tu,v,Tw)=\rho^r(u,Tv,Tw)=[u, v, w]_\mathfrak{h}
\end{align*}
for all $x,y,z\in \mathfrak{g}$ and $u,v,w\in \mathfrak{h}$.
Then, $T$ is a homomorphic embedding tensor.
\end{exam}

In \cite{Pei}, the authors introduced the concept of triLeibniz algebra within the context of studying tri-average operators on Lie algebras. Building on their groundwork, we propose the notion of 3-tri-Leibniz dialgebras.

\begin{definition}
A   3-tri-Leibniz dialgebra is a tuple $(\mathfrak{g}, [\cdot,\cdot,\cdot]_\mathfrak{g},[\cdot,\cdot, \cdot]_\vdash, [\cdot, \cdot,\cdot]_\dashv, [\cdot, \cdot, \cdot]_\perp)$ consisting of a 3-Leibniz algebra $(\mathfrak{g}, [\cdot,\cdot,\cdot]_\mathfrak{g})$
 and a  3-tri-Leibniz  algebra $(\mathfrak{g},  [\cdot, \cdot,\cdot]_\vdash, [\cdot,\cdot, \cdot]_\dashv, [\cdot,\cdot, \cdot]_\perp)$ such that
 for all $x,y,z,a,b\in \mathfrak{g}$ and $\triangle\in \{\vdash,\dashv,\perp\}$,
 \begin{align}
 &[a,b,[x,y,z]_\mathfrak{g}]_{\vdash}=[[a,b,x]_{\vdash},y,z]_\mathfrak{g}+[x,[a,b,y]_{\vdash},z]_\mathfrak{g}+[x,y,[a,b,z]_{\vdash}]_\mathfrak{g},\label{4.7}\\
 &[[a,b,x]_{\perp},y,z]_\mathfrak{g}+[x,[a,b,y]_{\perp},z]_\mathfrak{g}+[x,y,[a,b,z]_{\perp}]_\mathfrak{g}=0,\label{4.8}\\
 &[[a,b,x]_{\dashv},y,z]_\mathfrak{g}+[x,[a,b,y]_{\dashv},z]_\mathfrak{g}+[x,y,[a,b,z]_{\dashv}]_\mathfrak{g}=0,\label{4.9}\\
 &[a,b,[x,y,z]_\mathfrak{g}]_{\vdash}=[x,y,[a,b,z]_\vdash]_{\mathfrak{g}},~~~~[a,[x,y,z]_\mathfrak{g},b]_{\perp}=[x,y,[a,z,b]_\perp]_{\mathfrak{g}},\label{4.10}\\
 &[[x,y,z]_\mathfrak{g},a,b]_{\dashv}=[x,y,[z,a,b]_\dashv]_{\mathfrak{g}},~~~~ [[x,y,z]_\mathfrak{g},a,b]_{\vdash}= [[x,y,z]_{\triangle} ,a,b]_{\vdash},\label{4.11}\\
 &[a,[x,y,z]_\mathfrak{g},b]_{\vdash}= [a, [x,y,z]_{\triangle},,b]_{\vdash},~~ ~~[[x,y,z]_\mathfrak{g},a, b]_{\perp}= [[x,y,z]_{\triangle},a, b]_{\perp},\label{4.12}\\
  &[a, b, [x,y,z]_\mathfrak{g}]_{\perp}= [a, b, [x,y,z]_{\triangle}]_{\perp},~~~~[a,b,[x,y,z]_\mathfrak{g}]_{\dashv}= [a,b,[x,y,z]_{\triangle}]_{\dashv},\label{4.13}\\
    & [a,[x,y,z]_\mathfrak{g}, b]_{\dashv}= [a,[x,y,z]_{\triangle}, b]_{\dashv}.\label{4.14}
\end{align}
\end{definition}

\begin{theorem}
Let $T:\mathfrak{h}\rightarrow \mathfrak{g}$ be  a homomorphic embedding tensor on the
 3-Leibniz algebra  $(\mathfrak{g},[\cdot, \cdot, \cdot]_\mathfrak{g})$  with respect to  a   action $(\mathfrak{h},[\cdot, \cdot, \cdot]_\mathfrak{h};\rho^l,\rho^m,\rho^r)$.
Then $(\mathfrak{h},[\cdot, \cdot, \cdot]_\mathfrak{h},[\cdot,\cdot,\cdot]_{\vdash}^T,[\cdot,\cdot,\cdot]_{\perp}^T,$ $[\cdot,\cdot,\cdot]_{\dashv}^T)$ is a 3-tri-Leibniz dialgebra, where
\begin{align*}
[u,v,w]_{\vdash}^T=\rho^l(Tu,Tv,w),~~[u,v,w]_{\perp}^T=\rho^m(Tu,v,Tw)~\text{and}~[u,v,w]_{\dashv}^T=\rho^r(u,Tv,Tw)
\end{align*}
for all $u,v,w\in \mathfrak{h}$.
\end{theorem}

\begin{proof}
We have already $(\mathfrak{h},[\cdot, \cdot, \cdot]_\mathfrak{h})$ is a 3-Leibniz algebra. Moreover, it follows from Proposition \ref{prop:ET} that
$(\mathfrak{h}, [\cdot,\cdot,\cdot]_{\vdash}^T,[\cdot,\cdot,\cdot]_{\perp}^T,[\cdot,\cdot,\cdot]_{\dashv}^T)$ is a  3-tri-Leibniz algebra. Now, for any  $u,v,w,s,t\in \mathfrak{h}$, according to Eqs.  \eqref{4.1}-\eqref{4.6}, we have
\begin{align*}
&[s,t,[u,v,w]_\mathfrak{h}]_{\vdash}^T=\rho^l(Ts,Tt,[u,v,w]_\mathfrak{h})\\
&=[\rho^l(Ts,Tt,u),v,w]_\mathfrak{h}+[u,\rho^l(Ts,Tt,v),w]_\mathfrak{h}+[u,v,\rho^l(Ts,Tt,w)]_\mathfrak{h}\\
&=[[s,t,u]_{\vdash}^T,v,w]_\mathfrak{h}+[u,[s,t,v]_{\vdash}^T,w]_\mathfrak{h}+[u,v,[s,t,w]_{\vdash}^T]_\mathfrak{h},\\
&[[s,t,u]_{\perp}^T,v,w]_\mathfrak{h}+[u,[s,t,v]_{\perp}^T,w]_\mathfrak{h}+[u,v,[s,t,w]_{\perp}^T]_\mathfrak{h}\\
&=[\rho^m(Ts,t,Tu),v,w]_\mathfrak{h}+[u,\rho^m(Ts,t,Tv),w]_\mathfrak{h}+[u,v,\rho^m(Ts,t,Tw)]_\mathfrak{h}=0,\\
&[[s,t,u]_{\dashv}^T,v,w]_\mathfrak{h}+[u,[s,t,v]_{\dashv}^T,w]_\mathfrak{h}+[u,v,[s,t,w]_{\dashv}^T]_\mathfrak{h}\\
&=[\rho^r(s,Tt,Tu),v,w]_\mathfrak{h}+[u,\rho^r(s,Tt,Tv),w]_\mathfrak{h}+[u,v,\rho^r(s,Tt,Tw)]_\mathfrak{h}=0,\\
&[s,t,[u,v,w]_\mathfrak{h}]_{\vdash}^T=\rho^l(Ts,Tt,[u,v,w]_\mathfrak{h})=[u,v,\rho^l(Ts,Tt,w)]_\mathfrak{h}=[u,v,[s,t,w]_{\vdash}^T]_\mathfrak{h},\\
&[s,[u,v,w]_\mathfrak{h},t]_{\perp}^T=\rho^m(Ts,[u,v,w]_\mathfrak{h},Tt)=[u,v,\rho^m(Ts,w, Tt)]_\mathfrak{h}=[u,v,[s,w,t]_{\perp}^T]_\mathfrak{h},\\
&[[u,v,w]_\mathfrak{h},s,t]_{\dashv}^T=\rho^r([u,v,w]_\mathfrak{h},Ts,Tt)=[u,v,\rho^r(w,Ts, Tt)]_\mathfrak{h}=[u,v,[w, s,t]_{\dashv}^T]_\mathfrak{h}.
\end{align*}
On the other hand, we have
 \begin{align*}
&[[u,v,w]_\mathfrak{h},s,t]_{\vdash}^T =\rho^l(T[u,v,w]_\mathfrak{h},Ts,t)=\rho^l([Tu,Tv,Tw]_\mathfrak{g},Ts,t)=\rho^l(T[u,v,w]_{\triangle}^T,Ts,t)\\
&=[[u,v,w]_{\triangle}^T,s,t]_{\vdash}^T.
\end{align*}
Similarly, we obtain
 \begin{align*}
&[s,[u,v,w]_\mathfrak{h},t]_{\vdash}^T =[s,[u,v,w]_{\triangle}^T,t]_{\vdash}^T,~~ [s,[u,v,w]_\mathfrak{h},t]_{\dashv}^T =[s,[u,v,w]_{\triangle}^T,t]_{\dashv}^T,\\
& [s,t,[u,v,w]_\mathfrak{h}]_{\dashv}^T =[s,t,[u,v,w]_{\triangle}^T]_{\dashv}^T,~~~ [s,t,[u,v,w]_\mathfrak{h}]_{\perp}^T =[s,t,[u,v,w]_{\triangle}^T]_{\perp}^T,\\
&  [[u,v,w]_\mathfrak{h}s,t]_{\perp}^T =[[u,v,w]_{\triangle}^T, s,t]_{\perp}^T.
\end{align*}
Thus, Eqs  \eqref{4.7}-\eqref{4.14} are satisfied.
\end{proof}

\begin{coro}
Let $(\mathfrak{g},[\cdot, \cdot,\cdot]_\mathfrak{g})$ be a 3-Leibniz algebra and $T: \mathfrak{g}\rightarrow  \mathfrak{g}$  be a homomorphic averaging operator. Then,
$(\mathfrak{g},[\cdot, \cdot,\cdot]_\mathfrak{g},[\cdot,\cdot,\cdot]_{\vdash}^T,[\cdot,\cdot,\cdot]_{\perp}^T,[\cdot,\cdot,\cdot]_{\dashv}^T)$ is a     3-tri-Leibniz dialgebra,   where
 \begin{align*}
[x,y,z]_{\vdash}^T=[Tx,Ty,z]_\mathfrak{g},~~[x,y,z]_{\perp}^T=[Tx,y,Tz]_\mathfrak{g}~\text{and}~[x,y,z]_{\dashv}^T=[x,Ty,Tz]_\mathfrak{g}.
\end{align*}
\end{coro}

\section{  Linear deformations of embedding tensors on 3-Leibniz algebras }\label{sec: Deformations}
\def\theequation{\arabic{section}.\arabic{equation}}
\setcounter{equation} {0}

 This section  introduces 1-cocycles and the first cohomology  of embedding tensors on 3-Leibniz algebras.  Subsequently, we examine linear deformations of embedding tensors on 3-Leibniz algebras,
  adhering to Gerstenhaber's methodology \cite{Gerstenhaber}.
In particular, we introduce the notion of a Nijenhuis element associated with an   embedding tensor, which leads to a trivial linear deformations.

Let $T:V\rightarrow  \mathfrak{g}$ be   an  embedding tensor  on a 3-Leibniz algebra  $(\mathfrak{g},[\cdot,\cdot,\cdot]_{\mathfrak{g}})$  with respect to the representation $(V; \rho^l,\rho^m,\rho^r)$.

 \begin{definition}
A parameterized sum $T_t=T+tT_1$, for some $T_1\in \mathrm{Hom}(V, \mathfrak{g})$, is called a  linear deformation of $T$ if $T_t$ is an   embedding tensor
 for all values of parameter $t$. In this case, we say that $T_1$ generates a linear deformation of $T$.
\end{definition}

Suppose that $T_1$  generates a linear deformation of  $T$, then we have
$$ [T_tu,T_tv,T_tw]_\mathfrak{g}=T_t\rho^l(T_tu,T_tv,w)=T_t\rho^m(T_tu,v,T_tw)=T_t\rho^r(u,T_tv,T_tw),$$
for all $u,v,w\in V$. This is equivalent to the following equations
 \begin{align}
&[Tu,Tv,T_1w]_\mathfrak{g}+[Tu,T_1v,Tw]_\mathfrak{g}+[T_1u,Tv,Tw]_\mathfrak{g}\nonumber\\
&=
\left\{ \begin{array}{lll}
=T_1\rho^l(Tu,Tv,w)+T\rho^l(T_1u,Tv,w)+T\rho^l(Tu,T_1v,w),\\
=T_1\rho^m(Tu,v,Tw)+T\rho^m(T_1u,v,Tw)+T\rho^m(Tu,v,T_1w),\\
=T_1\rho^r(u,Tv,Tw)+T\rho^r(u,Tv,T_1w)+T\rho^r(u,T_1v,Tw),
 \end{array}  \right.\label{5.1}\\
&[Tu,T_1v,T_1w]_\mathfrak{g}+[T_1u,Tv,T_1w]_\mathfrak{g}+[T_1u,T_1v,Tw]_\mathfrak{g}\nonumber\\
&=\left\{ \begin{array}{lll}
=T\rho^l(T_1u,T_1 v,w)+T_1\rho^l(Tu,T_1v,w)+T_1\rho^l(T_1u,Tv,w),\\
=T\rho^m(T_1 u,v,T_1w)+T_1\rho^m(Tu,v,T_1w)+T_1\rho^m(T_1u,v,Tw),\\
=T\rho^r(u,T_1 v,T_1w)+T_1\rho^r(u,Tv,T_1w)+T_1\rho^r(u,T_1v,Tw),
 \end{array}  \right.\label{5.2}\\
&[T_1u,T_1v,T_1w]_\mathfrak{g}=T_1\rho^l(T_1u,T_1v,w)=T_1\rho^m(T_1u,v,T_1w)=T_1\rho^r(u,T_1v,T_1w).\label{5.3}
\end{align}
Thus, $T_t$ is a linear deformation of $T$ if and only if Eqs. \eqref{5.1}-\eqref{5.3} hold.
From Eq. \eqref{5.3} it follows that the map $T_1$ is an  embedding tensor  on a 3-Leibniz algebra  $(\mathfrak{g},[\cdot,\cdot,\cdot]_{\mathfrak{g}})$  with respect to the representation $(V; \rho^l,\rho^m,\rho^r)$.

\begin{definition}
A linear map  $T_1\in \mathrm{Hom}(V, \mathfrak{g})$  is called a 1-cocycle of the   embedding tensor $T$  if it satisfies Eq. \eqref{5.1}.
One denotes by $\mathbf{Z}_{T}^1(V, \mathfrak{g})$ the set of theses 1-cocycles.

Moreover, two  1-cocycles $T_1$ and $\widetilde{T}_1$ are said to be equivalent
  if there
exists  an element $(a, b)\in \wedge^2 \mathfrak{g}$ such that for all  $u\in V$,
   \begin{align}
\widetilde{T}_1u-T_1u=T\rho^l(a,b,u)-[a,b,Tu]_\mathfrak{g}=\delta(a,b)u. \label{5.4}
\end{align}
In this case, we denote $T_1\sim \widetilde{T}_1$.
The first    cohomology   $\mathrm{H}\mathrm{H}_{T}^1(V, \mathfrak{g})$ is the quotient of $\mathbf{Z}_{T}^1(V, \mathfrak{g})$  by this equivalence relation.
\end{definition}

Let  $(\mathfrak{g},[\cdot, \cdot, \cdot]_\vdash, [\cdot, \cdot, \cdot]_\dashv, [\cdot, \cdot, \cdot]_\perp)$   be a 3-tri-Leibniz algebra,
 and let $\omega_\vdash, \omega_\dashv$ and $\omega_\perp$ be   trilinear maps on $\mathfrak{g}$. If for any values of parameter $t$,
the brackets $[\cdot, \cdot, \cdot]^t_\vdash, [\cdot, \cdot, \cdot]^t_\dashv, [\cdot, \cdot, \cdot]^t_\perp$ defined by
 \begin{align*}
 &[x,y,z]^t_{\vdash}=[x,y,z]_\vdash+t\omega_\vdash(x,y,z),\\
 &[x,y,z]^t_{\dashv}=[x,y,z]_\dashv+t\omega_\dashv(x,y,z),\\
 &[x,y,z]^t_{\perp}=[x,y,z]_\perp+t\omega_\perp(x,y,z)
  \end{align*}
for all $x,y,z\in \mathfrak{g}$,  also gives a  3-tri-Leibniz algebra structure, we say that the
tuple $(\omega_\vdash, \omega_\dashv,\omega_\perp)$ generates a linear deformation of
the 3-tri-Leibniz algebra $(\mathfrak{g},[\cdot, \cdot, \cdot]_\vdash, [\cdot, \cdot, \cdot]_\dashv, [\cdot, \cdot, \cdot]_\perp)$. The two types of linear deformations are related as follows.

  \begin{prop}
If $T_1$  generates a linear deformation of  $T$, then  the
tuple $(\omega_\vdash, \omega_\dashv,\omega_\perp)$   defined by
 \begin{align*}
 &\omega_{\vdash}(u,v,w)=\rho^l(T_1u,Tv,w)+\rho^l(Tu,T_1v,w),\\
 &\omega_{\perp}(u,v,w)=\rho^m(T_1u,v,Tw)+\rho^m(Tu,v,T_1w),\\
 &\omega_{\dashv}(u,v,w)=\rho^r(u,Tv,T_1w)+\rho^r(u,T_1v,Tw)
  \end{align*}
for all $u,v,w\in  V$, generates a linear deformation of the  3-tri-Leibniz algebra $(V,[\cdot,\cdot,\cdot]_{\vdash}^T,$ $[\cdot,\cdot,\cdot]_{\perp}^T,[\cdot,\cdot,\cdot]_{\dashv}^T)$ given  in
Proposition \ref{prop:ET}.
\end{prop}

\begin{proof}
Let $T_t=T+tT_1$ be a linear deformation of $T$.
According to Proposition \ref{prop:ET}, we denote by $([\cdot, \cdot, \cdot]^{T_t}_\vdash, [\cdot, \cdot, \cdot]^{T_t}_\dashv, [\cdot, \cdot, \cdot]^{T_t}_\perp)$
the 3-tri-Leibniz algebra structure on $V$ corresponding  to the  embedding tensor $T_t$.
Then, for any $u,v,w\in  V$, we have
  \begin{align*}
&[u, v, w]^{T_t}_\vdash= \rho^l(T_tu,T_tv,w)\\
&= \rho^l(Tu,Tv,w)+ \rho^l(T_1u,Tv,w)+\rho^l(Tu,T_1v,w)~~~( \mathrm{mod}~t^2)\\
&=[u,v,w]^T_\vdash+t\omega_\vdash(u,v,w).
\end{align*}
Similarly, we obtain
  \begin{align*}
&[u, v, w]^{T_t}_\perp=[u,v,w]^T_\perp+t\omega_\perp(u,v,w),\\
&[u, v, w]^{T_t}_\dashv=[u,v,w]^T_\dashv+t\omega_\dashv(u,v,w),
\end{align*}
which implies that  $(\omega_\vdash, \omega_\dashv,\omega_\perp)$  generates a linear deformation of  $(V,[\cdot,\cdot,\cdot]_{\vdash}^T,[\cdot,\cdot,\cdot]_{\perp}^T,$ $[\cdot,\cdot,\cdot]_{\dashv}^T)$.
\end{proof}

\begin{definition}
Let $T$ and $\widetilde{T}$
be two embedding tensors on the 3-Leibniz algebra  $(\mathfrak{g},[\cdot,\cdot,\cdot]_{\mathfrak{g}})$  with respect to the representation $(V; \rho^l,\rho^m,\rho^r)$.
 A homomorphism from $\widetilde{T}$
to $T$ is a pair $(\varphi, \phi)$, where  $\varphi:\mathfrak{g}\rightarrow\mathfrak{g}$ is a 3-Leibniz algebra homomorphism and $\phi:V\rightarrow V$
is a linear map such that for all $x,y\in\mathfrak{g},  u\in V$,
  \begin{align}
&\phi(\rho^l(x,y,u))=\rho^l(\varphi(x),\varphi(y),\phi(u)),\label{5.5}\\
&\phi(\rho^m(x,u,y))=\rho^m(\varphi(x),\phi(u),\varphi(y)),\label{5.6}\\
&\phi(\rho^r(u,x,y))=\rho^m(\phi(u),\varphi(x),\varphi(y)),\label{5.7}\\
& T \phi(u)= \varphi(\widetilde{T}u).\label{5.8}
\end{align}
\end{definition}

 \begin{definition}
Two linear deformations $T_t=T+tT_1$ and $\widetilde{T}_t=T+t\widetilde{T}_1$
of an  embedding tensor $T$  on a 3-Leibniz algebra  $(\mathfrak{g},[\cdot,\cdot,\cdot]_{\mathfrak{g}})$  with respect to the representation $(V; \rho^l,\rho^m,\rho^r)$  are said to
be equivalent if there exists an element $\mathfrak{A}=(a, b)\in \wedge^2 \mathfrak{g}$ such that the pair
$$(\varphi_t=\mathrm{id}_{ \mathfrak{g}}+t[a,b,\cdot]_\mathfrak{g},\phi_t=\mathrm{id}_{V}+t\rho^l(a,b,\cdot))$$
defines a homomorphism of  embedding tensors
from $\widetilde{T}_t$ to $T_t$.

In particular, a linear deformation $T_t=T+tT_1$ of an embedding
tensor $T$ is called trivial if there exists an element  $(a, b)\in \wedge^2 \mathfrak{g}$ such that $(\varphi_t,\phi_t)$ is a
homomorphism from $T$ to $T_t$.
 \end{definition}

Let $(\varphi_t,\phi_t)$  be a homomorphism from from $\widetilde{T}_t$ to $T_t$.  Then $\varphi_t$ is  a  3-Leibniz algebra   homomorphism, which implies that
\begin{align}
\left\{ \begin{array}{ll}
~[[a, b, x]_\mathfrak{g},[a, b, y]_\mathfrak{g},z]_\mathfrak{g}+ [[a, b, x]_\mathfrak{g},y, [a, b, z]_\mathfrak{g}]_\mathfrak{g}+ [x,[a, b, y]_\mathfrak{g},[a, b, z]_\mathfrak{g}]_\mathfrak{g}=0,\\
~[[a, b, x]_\mathfrak{g},[a, b, y]_\mathfrak{g},[a, b, z]_\mathfrak{g}]_\mathfrak{g}=0
 \end{array}  \right. \label{5.9}
 \end{align}
for all $x,y,z\in  \mathfrak{g}$.

Note that by Eqs. \eqref{5.5}--\eqref{5.7}, we obtain that for all $x,y\in\mathfrak{g},  u\in V,$
\begin{small}
\begin{align}
&\left\{ \begin{array}{lll}
\rho^l\big(a,b,\rho^l(x,y,u)\big)=\rho^l\big([a,b,x]_\mathfrak{g},y,u \big)+\rho^l\big(x,[a,b,y]_\mathfrak{g},u \big)+\rho^l\big(x,y,\rho^l(a,b,u)\big),\\
\rho^l\big(a,b,\rho^m(x,u,y)\big)=\rho^m\big([a,b,x]_\mathfrak{g},u,y \big)+\rho^m\big(x,\rho^l(a,b,u),y \big)+\rho^m\big(x,u,[a,b,y]_\mathfrak{g}\big),\\
\rho^l\big(a,b,\rho^r(u,x,y)\big)=\rho^r\big(\rho^l(a,b,u),x,y \big)+\rho^r\big(u,[a,b,x]_\mathfrak{g},y \big)+\rho^m\big(u,x,[a,b,y]_\mathfrak{g}\big),
 \end{array}  \right. \label{5.10}\\
&\left\{ \begin{array}{lll}
\rho^l(x,[a,b,y]_\mathfrak{g},\rho^l(a,b,u) )+\rho^l([a,b,x]_\mathfrak{g},y,\rho^l(a,b,u))+\rho^l([a,b,x]_\mathfrak{g},[a,b,y]_\mathfrak{g},u )=0,\\
\rho^m(x,\rho^l(a,b,u),[a,b,y]_\mathfrak{g})+\rho^m([a,b,x]_\mathfrak{g},u,[a,b,y]_\mathfrak{g} )+\rho^m([a,b,x]_\mathfrak{g},\rho^l(a,b,u),y )=0,\\
\rho^r(u,[a,b,x]_\mathfrak{g},[a,b,y]_\mathfrak{g} )+\rho^m(\rho^l(a,b,u),x,[a,b,y]_\mathfrak{g})+\rho^r(\rho^l(a,b,u),[a,b,x]_\mathfrak{g},y )=0,
 \end{array}  \right. \label{5.11}\\
&\left\{ \begin{array}{lll}
\rho^l([a,b,x]_\mathfrak{g},[a,b,y]_\mathfrak{g},\rho^l(a,b,u) )=0,\\
\rho^m([a,b,x]_\mathfrak{g},\rho^l(a,b,u),[a,b,y]_\mathfrak{g})=0,\\
\rho^r(\rho^l(a,b,u),[a,b,x]_\mathfrak{g},[a,b,y]_\mathfrak{g} )=0.
 \end{array}  \right. \label{5.12}
 \end{align}
 \end{small}

Moreover, Eq.   \eqref{5.8}   yields that for all  $u\in V$,
 \begin{align*}
&T_t(u+t\rho^l(a,b,u))= \varphi_t(Tu+t\widetilde{T}_1u),
\end{align*}
which implies
  \begin{align}
T_1u+T\rho^l(a,b,u)=&  \widetilde{T}_1u+[a,b,Tu]_\mathfrak{g}, \label{5.13}\\
T_1\rho^l(a,b,u)=&[a,b,\widetilde{T}_1u]_\mathfrak{g}.\nonumber
\end{align}

This shows that $T_1$ and  $\widetilde{T}_1$ are cohomologous. Hence their cohomology classes are the same in $\mathrm{H}\mathrm{H}_{T}^1(V, \mathfrak{g})$.
Conversely, any 1-cocycle $T_1$  gives rise to the linear deformation $T+tT_1$. Moreover, cohomologous
1-cocycles correspond to equivalent linear deformations.
In summary of the aforementioned discussions, we arrive at the following conclusion.

\begin{theorem}
Let $T:V\rightarrow  \mathfrak{g}$ be   an  embedding tensor on a 3-Leibniz algebra  $(\mathfrak{g},[\cdot,\cdot,\cdot]_{\mathfrak{g}})$  with respect to the representation $(V; \rho^l,\rho^m,\rho^r)$.
 Then there is a bijection between the set of all
equivalence classes of linear deformations of $T$ and the first cohomology   $\mathrm{H}\mathrm{H}_{T}^1(V, \mathfrak{g})$.
\end{theorem}

 \begin{definition}
Let $T:V\rightarrow  \mathfrak{g}$ be   an  embedding tensor on a 3-Leibniz algebra  $(\mathfrak{g},[\cdot,\cdot,\cdot]_{\mathfrak{g}})$  with respect to the representation $(V; \rho^l,\rho^m,\rho^r)$.
 An element $\mathfrak{A}=(a, b)\in \wedge^2 \mathfrak{g}$ is said to be a Nijenhuis element associated to $T$ if $\mathfrak{A}$
satisfies Eqs. \eqref{5.9}-\eqref{5.12} and the equation
  \begin{align}
&[a,b,T\rho^l(a,b,u)-[a,b,Tu]_\mathfrak{g}]_\mathfrak{g}=0  \label{5.14}
\end{align}
for all $u\in  V$.  Denote by $\mathrm{Nij}(T)$ the set of Nijenhuis elements associated with an embedding tensor $T$.
 \end{definition}

It is easy to see that a trivial linear deformation of an  embedding tensor   on a 3-Leibniz algebra  gives rise to a Nijenhuis element. Nonetheless, the reverse holds equally true.

 \begin{theorem} \label{theorem:TID}
Let $T:V\rightarrow  \mathfrak{g}$ be   an  embedding tensor on a 3-Leibniz algebra  $(\mathfrak{g},[\cdot,\cdot,\cdot]_{\mathfrak{g}})$  with respect to the representation $(V; \rho^l,\rho^m,\rho^r)$.
 Then, for any $\mathfrak{A}=(a, b)\in \mathrm{Nij}(T)$, $T_t=T+tT_1$ with $T_1=\delta(a,b)$ is
a trivial  linear deformation of the   embedding tensor $T$.
\end{theorem}

 To establish the proof for this theorem, we require the subsequent lemma.

\begin{lemma} \label{lemma:newET}
Let $T:V\rightarrow  \mathfrak{g}$ be   an  embedding tensor on a 3-Leibniz algebra  $(\mathfrak{g},[\cdot,\cdot,\cdot]_{\mathfrak{g}})$  with respect to the representation $(V; \rho^l,\rho^m,\rho^r)$.
  Let $\varphi: \mathfrak{g}\rightarrow\mathfrak{g} $ be a 3-Leibniz algebra  isomorphism   and $\phi:V\rightarrow V $ be a reversible linear map   such that  Eqs. \eqref{5.5}--\eqref{5.7} hold.
  Then $\varphi^{-1}\circ T\circ \phi:V\rightarrow  \mathfrak{g}$ is an embedding tensor
   on a 3-Leibniz algebra  $(\mathfrak{g},[\cdot,\cdot,\cdot]_{\mathfrak{g}})$  with respect to the representation $(V; \rho^l,\rho^m,\rho^r)$.
\end{lemma}

\begin{proof}
For any $u,v,w\in V$,  we have
\begin{align*}
&[(\varphi^{-1}\circ T\circ \phi)(u),(\varphi^{-1}\circ T\circ \phi)(v),(\varphi^{-1}\circ T\circ \phi)(w)]_\mathfrak{g}=\varphi^{-1}[T\phi(u), T\phi(v),  T \phi(w)]_\mathfrak{g}\\
&=\left\{ \begin{array}{lll}
 =\varphi^{-1}T\rho^l(T\phi(u), T\phi(v),  \phi(w))=(\varphi^{-1}\circ T\circ \phi)\rho^l(\varphi^{-1}(T\phi(u)), \varphi^{-1}(T\phi(v)),   w),\\
 =\varphi^{-1}T\rho^m(T\phi(u), \phi(v),  T\phi(w)) =(\varphi^{-1}\circ T\circ \phi)\rho^m(\varphi^{-1}(T\phi(u)), v,  \varphi^{-1}(T\phi(w))),\\
 =\varphi^{-1}T\rho^r(\phi(u), T\phi(v),  T\phi(w))=(\varphi^{-1}\circ T\circ \phi)\rho^r(u, \varphi^{-1}(T\phi(v)),   \varphi^{-1}(T\phi(w))),
 \end{array}  \right.
\end{align*}
which implies that   $\varphi^{-1}\circ T\circ \phi$ is an embedding tensor.
\end{proof}

The proof of Theorem \ref{theorem:TID}:
 For any Nijenhuis element $\mathfrak{A}=(a, b)\in \mathrm{Nij}(T)$, we define a linear map $T_1:V\rightarrow \mathfrak{g}$ by
$$T_1u=\delta(a,b)u= T\rho^l(a,b,u)-[a,b,Tu]_\mathfrak{g} $$
for all $u\in V$.
Let $T_t=T+tT_1$. By the Eqs. \eqref{5.9}-\eqref{5.12} and \eqref{5.14}, for all $x,y,z\in  \mathfrak{g}$ and  $u\in V$, we have
 \begin{align*}
[x,y,z]_{\mathfrak{g}}+t[a,b,[x,y,z]_{\mathfrak{g}}]_{\mathfrak{g}}&=[x+t[a,b,x]_{\mathfrak{g}},y+t[a,b,y]_{\mathfrak{g}},z+t[a,b,z]_{\mathfrak{g}}]_{\mathfrak{g}},\\
\rho^l(x,y,u)+t\rho^l(a,b,\rho^l(x,y,u))&=\rho^l(x+t[a,b,x]_{\mathfrak{g}},y+t[a,b,y]_{\mathfrak{g}},u+t\rho^l(a,b,u)),\\
\rho^m(x,u,y)+t\rho^l(a,b,\rho^m(x,u,y))&=\rho^m(x+t[a,b,x]_{\mathfrak{g}},u+t\rho^l(a,b,u),y+t[a,b,y]_{\mathfrak{g}}),\\
\rho^r(u,x,y)+t\rho^l(a,b,\rho^r(u,x,y))&=\rho^r(u+t\rho^l(a,b,u),x+t[a,b,x]_{\mathfrak{g}},y+t[a,b,y]_{\mathfrak{g}}),\\
T_tu+t[a,b,T_tu]_\mathfrak{g}&=Tu+tT\rho^l(a,b,u).
\end{align*}
Since $[a,b,\cdot]_{\mathfrak{g}}$ and $\rho^l(a,b,\cdot)$ are linear transformations of finite-dimensional  vector spaces $ \mathfrak{g}$ and $V$, respectively, it follows that
$\varphi_t=\mathrm{id}_{ \mathfrak{g}}+t[a,b,\cdot]_\mathfrak{g}$ is a 3-Leibniz algebra  isomorphism, and $\phi_t=\mathrm{id}_{V}+t\rho^l(a,b,\cdot)$ is reversible linear map, respectively, for $|t|$  sufficiently small.
Moreover by  Eq. \eqref{5.14}, for all $u\in V$,  we have,
 \begin{align*}
&\varphi_t^{-1}\circ T\circ\phi_t(u)=(\mathrm{id}_{\mathfrak{g}}+t[a,b,\cdot]_\mathfrak{g})^{-1}\circ T\circ(\mathrm{id}_{V}+t\rho^l(a,b,\cdot))(u)\\
&=\sum_{i=0}^{+\infty}(-t[a,b,\cdot]_\mathfrak{g})^i\circ(Tu+tT\rho^l(a,b,u))\\
&=Tu+t(T\rho^l(a,b,u)-[a,b,Tu]_\mathfrak{g})+\sum_{i=1}^{+\infty}(-1)^it^{i+1}[a,b,\cdot]_\mathfrak{g}^i\circ(T\rho^l(a,b,u)-[a,b,Tu]_\mathfrak{g})\\
&=Tu+t(T\rho^l(a,b,u)-[a,b,Tu]_\mathfrak{g})\\
&=Tu+tT_1u=T_tu.
\end{align*}
So  by  Lemma \ref{lemma:newET}, $T_t$ is an  embedding tensor on a 3-Leibniz algebra  $(\mathfrak{g},[\cdot,\cdot,\cdot]_{\mathfrak{g}})$  with respect to the representation $(V; \rho^l,\rho^m,\rho^r)$ for $|t|$  sufficiently small.
Thus, $T_1$   satisfies the Eqs.  \eqref{5.1}-\eqref{5.3}. Therefore, $T_t$
is an  embedding tensor for all $t$, which means that $T_1$  generates
a linear deformation of $T$. It is straightforward to show that this linear deformation is
trivial.


\begin{thebibliography}{99}

\bibitem{Casas}  J. M. Casas, J. L. Loday,  T. Pirashvili.  Leibniz $n$-algebras. Forum Math., 2002, 14:  189--207.


\bibitem{Filippov} V. T. Filippov.  $n$-Lie algebras.  Sibirsk. Mat. Zh., 1985,  26:  126-140, 191.

 \bibitem{Albeverio} S.  Albeverio, S.  Ayupov, B.   Omirov, R.Turdibaev.  Cartan subalgebras of Leibniz $n$-algebras. Comm. Algebra, 2009, 37(6): 2080--2096.

\bibitem{Casas16}  J. M. Casas,  M.A. Insua, M.  Ladra,  S.  Ladra.  Test for Leibniz $n$-algebra structure. Linear Algebra Appl., 2016, 494:  138--155.

 \bibitem{XU}  N.  Xu, Y. Sheng. On non-abelian extensions of 3-Leibniz algebras. Front. Math. China,  2024, 19(2): 57--74.

 \bibitem{Azcarraga} J.A. de Azc$\acute{\mathrm{a}}$rraga,  J.M.  Izquierdo.  $n$-ary algebras: a review with applications. J. Phys. A , 2010, 43: 293001.

 \bibitem{Nicolai}  H. Nicolai, H. Samtleben.  Maximal gauged supergravity in three dimensions.  Phys. Rev. Lett., 2001,  86:  1686--1689.

 \bibitem{Aguiar} M. Aguiar.   Pre-Poisson Algebras.  Lett. Math. Phys., 2000,  54:  263--277.



\bibitem{Sheng}  Y. Sheng,  R.Tang,  C.Zhu. The controlling $L_\infty$-algebra, cohomology and homotopy of embedding tensors and Lie-Leibniz triples. Commun. Math. Phys., 2021, 386:  269--304.

\bibitem{Braiek}  S. Braiek, T. Chtioui,   S. Mabrouk. Anti-Leibniz algebras: A non-commutative version of mock-Lie algebras.	arXiv:2409.17184.

\bibitem{Das24}A. Das.   Averaging operators on groups, racks and Leibniz algebras. 	arXiv:2403.06250.

\bibitem{Caseiro} R.  Caseiro, J. M. Costa.  Embedding tensors on Lie$_\infty$-algebras with respect to Lie$_\infty$-actions.  Commun. Algebra, 2024,  52(4):  1432--1456.

 \bibitem{Hu}  M. Hu, S. Hou, L. Song, et al.   Deformations and cohomologies of embedding tensors on 3-Lie algebras.  Comm. Algebra,  2023, 51(11):  4622--4639.

\bibitem{Teng25} W. Teng.  Embedding tensors on Lie triple systems. To appear in  Filomat. 


\bibitem{Gerstenhaber}  M. Gerstenhaber. On the deformation of rings and algebras. Ann. Math.,  1964, 2(79):  59--103.


\bibitem{Loday} J.  Loday. Dialgebras, in ``Dialgebras and related operads". Springer Lecture Notes in Math., 2001, 1763: 7--66.



\bibitem{Pei}  J. Pei,  C. Bai, L. Guo,  X. Ni. Replicators, Manin white product of binary operads and average operators. Conference: 3rd International Congress in Algebras and Combinatorics, 2020.

































































































\end{thebibliography}
\end{document}